\documentclass[12pt]{article}
\usepackage[pctex32]{graphics}
\usepackage{amsthm,amsmath,amssymb,amscd,verbatim,epsfig}
\usepackage{graphics,pifont}
\usepackage{amssymb}         %
\usepackage{graphicx}
\newcommand{\beq}{\begin{equation}}
\newcommand{\eeq}{\end{equation}}

\date{}

\newcommand{\supp}{\supset}

\newcommand{\sq}{$\blacksquare$}

\begin{document}

\title{A\ Simple\ Proof\ of\ Sharkovsky's\ Theorem\ Revisited}
\author{Bau-Sen Du \\ [.3cm]
Institute of Mathematics \\
Academia Sinica \\
Taipei 11529, Taiwan \\
dubs@math.sinica.edu.tw \\}
\maketitle


Sharkovsky's theorem {\bf{\cite{bl}}}, {\bf{\cite{mi}}}, {\bf{\cite{sh}}} claims its place among the gems of dynamical systems mainly because of its simple hypotheses and strong conclusion.  Loosely speaking, it states that, if $f$ is a continuous map from a compact interval $I$ into itself that has a period-$m$ point, then $f$ also has a period-$n$ point whenever $m \prec n$ in the Sharkovsky's ordering $\prec$ of the natural numbers but may not have a period-$k$ point when $k \prec m$.  The Sharkovsky's ordering is as follows: $$3 \prec 5 \prec 7 \prec  9 \prec \cdots \prec 2 \cdot 3 \prec 2 \cdot 5 \prec 2 \cdot 7 \prec  2 \cdot 9 \prec \cdots \prec 2^2 \cdot 3 \prec 2^2 \cdot 5 \prec 2^2 \cdot 7 \prec  2^2 \cdot 9 \prec \cdots$$ $$\prec \cdots \prec 2^3 \prec 2^2 \prec 2 \prec 1.$$  In {\bf{\cite{du}}}, a simple proof of Sharkovsky's theorem using directed graphs was presented.  In the present note, we give another simple proof of the theorem without using directed graphs. 

\begin{figure}[htb]
\centerline{\epsfig{file=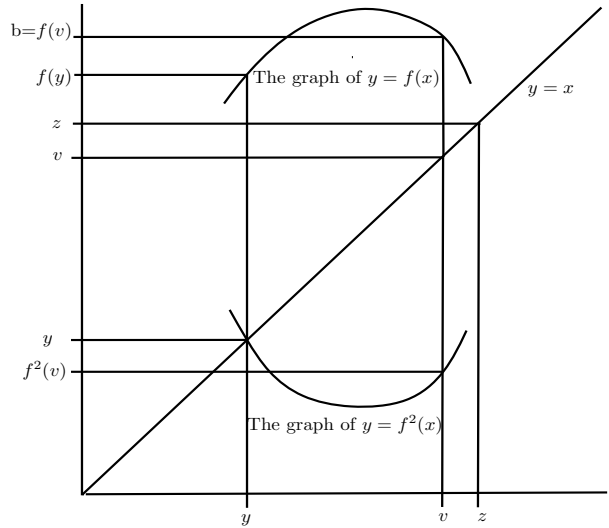,width=8cm,height=7cm}}
\caption{The relative locations of the points $v$, $f(v)$, $f^2(v)$, $y$, $f(y)$, and $z$.}
\end{figure}

We need the following elementary result, which may at first glance look a little complicated.  However, it can be best illustrated by Figure 1.  With the help of Figure 1, it is also easier to see why Theorems 2 and 3 follow from Lemma 1 so naturally.

\noindent
{\bf Lemma 1.}
{\it Let $a$ and $b$ be points of $I$ such that either $f(b) < a < b \le f(a)$ or $f(b) \le a < b < f(a)$.  Then there exist a fixed point $z$ of $f$ with $z < b$, a period-2 point $y$ of $f$ with $y < z$, and a point $v$ in $(y, z)$ such that $f(v) = b$ and $$\max \{ f^2(v), y \} < v < z < \min \{ f(y), f(v) \}.$$  Furthermore, $f(x) > z$ and $f^2(x) < x$ whenever $y < x \le v$.}

\noindent
{\it Proof.}
No matter whether $f(b) < a < b \le f(a)$ or $f(b) \le a < b < f(a)$, $f$ has a fixed point $z$ in $(a, b)$ and a point $v$ in $[a, z)$ such that $f(v) = b > z$.  If $f(x) > z$ whenever $\min I \le x \le v$, let $u = \min I$; otherwise let $u = \max \{x : \min I \le x \le v$, $f(x) = z \}$.  Then $f^2(u) \ge u$ and $f(x) > z$ when $u < x \le v$.  Since $f^2(v) < v$, $f$ has a period-2 point in $[u, v)$.  If $y$ is the largest such point, then $u \le y < v < z < f(y)$.  Finally, since $f^2(v) < v$, $f^2(x) < x$ for each $x$ in $(y, v]$.  
\hfill\sq

Let $P$ be a period-$m$ orbit of $f$ with $m \ge 3$, let $p$ and $b$ $(p < b)$ be points in $P$ such that $f(p) \ge b$ and $f(b)\le p$ (so $f$ has a fixed point in $[p, b]$), and let $a$ be a point in $[p, b)$ such that $f(a) = b$.  Since $f(b) < a \, (< b = f(a))$, we see that Lemma 1 applies and note that $b$, as a point in $P$, has least period $m$.  These facts will be used in the proofs of Theorems 2 and 3.

\noindent
{\bf Theorem 2.}
{\it If $f$ has a periodic point of least period $m$ with $m \ge 3$ and odd, then $f$ has periodic points of least period $n$ for each odd integer $n$ such that $n \ge m$.}

\noindent
{\it Proof.}
Let $P$ be a period-$m$ orbit of $f$.  By Lemma 1 there exist a fixed point $z$ of $f$, a period-2 point $y$ of $f$ and a point $v$ with $y < v < z < f(y)$ such that $f(v)$ lies in $P$ and, $f(x) > z$ and $f^2(x) < x$ when $y < x \le v$.  Let $p_m = v$.  Because $m$ is odd, $f^{m+2}(y) = f(y) > y$, and because $f^2(p_m) \, (= f^2(v))$ is a period-$m$ point of $f$, $f^{m+2}(p_m) = f^2(p_m) < p_m$.  Thus the point $p_{m+2} = \min \{x : y \le x \le p_m$, $f^{m+2}(x) = x \}$ exists and is a period-$(m+2)$ point of $f$.  Since $f^{m+4}(y) = f(y) > y$ and $f^{m+4}(p_{m+2}) = f^2(p_{m+2}) < p_{m+2}$, the point $p_{m+4} = \min \{x : y \le x \le p_{m+2}$, $f^{m+4}(x) = x \}$ exists and is a period-$(m+4)$ point of $f$.  Proceeding in this manner indefinitely, we obtain a decreasing sequence $$y < \cdots < p_{m+2i+2} < p_{m+2i} < \cdots < p_{m+2} < p_m = v$$ of points such that $p_{m+2i}$ is a period-$(m+2i)$ point of $f$ ($i = 1, 2, \ldots$).   
\hfill\sq

\noindent
{\bf Theorem 3.}
{\it If $f$ has a periodic point of least period $m$ with $m \ge 3$ and odd, then $f$ has periodic points of all even periods and $f^2$ is strictly turbulent (i.e., there exist disjoint closed subintervals $I_0$ and $I_1$ of $I$ such that $f^2(I_0) \cap f^2(I_1) \supp I_0 \cup I_1$).}

\noindent
{\it Proof.}
Consider a period-$m$ orbit $P$ of $f$.  By Lemma 1, there exist a fixed point $z$ of $f$, a period-2 point $y$ of $f$, and a point $v$ such that $f(v) = b \in P$, $$\max\{f^2(v), y\} < v < z < b = f(v) = f^{m+1}(v),$$ and $f^2(x) < x$ and $f(x) > z$ for all $x$ in $(y, v]$.  Let $g = f^2$ and let $z_0 = \min \{ x: v \le x \le z$, $g(x) = x \}$.  Then $y$ and $z_0$ are fixed points of $g$ such that $y < v < z_0 \le z < b = g^{(m+1)/2}(v)$.  Furthermore, $g(x) < x$ and $f(x) > z$ when $y < x < z_0$.  If $g(x) < z_0$ whenever $\min I \le x < z_0$, then $g([\min I, z_0]) \subset [\min I, z_0]$, which contradicts the fact that $g^{(m+1)/2}(v) = b > z_0$.  Accordingly, the point $d = \max \{ x : \min I \le x \le y$, $g(x) = z_0 \}$ exists and $f(x) > z \ge z_0 > g(x)$ for all $x$ in $(d, y)$.  Therefore $f(x) > z \ge z_0 > g(x)$ whenever $d < x < z_0$.  Let $s = \min \{ g(x) : d \le x \le z_0 \}$.  If $s \ge d$, then $g([d, z_0]) \subset [d, z_0]$, which again contradicts the fact that $g^{(m+1)/2}(v) = b > z_0$.  Thus $s < d$.  Consequently, $g$ is strictly turbulent.  Let $\hat g(x)$ be the continous map from $[d, z_0]$ onto itself defined by $\hat g(x) = \max \{ g(x), d \}$, and let $t = \min \{ x: d \le x \le z_0$, $g(x) = d \}$.  For each integer $n \ge 1$ let $c_n = \min \{ x : d \le x \le t$, ${\hat g}^n(x) = x \}$.  It is easy to see that $d < \cdots < c_4 < c_3 < c_2 < c_1 \le y$ and that each $c_n$ generates a period-$n$ orbit $Q_n \subset (d, z_0)$ of $\hat g$.  It is obvious that each $Q_n$ is also a period-$n$ orbit of $g = f^2$.  Since $x < z_0 \le z < f(x)$ for every $x$ in $Q_n$, $Q_n \cup f(Q_n)$ is a period-$(2 \cdot n)$ orbit of $f$.  This shows that $f$ has periodic points of all even periods.
\hfill\sq

For the sake of completeness, we include a proof of Sharkovsky's theorem {\bf{\cite{mi}}} that is slightly different from the one presented in {\bf{\cite{du}}}.

\noindent
{\bf Theorem 4 (Sharkovsky).}
{\it The following three statements hold:
\begin{itemize}
\item[\rm{(1)}] 
Assume that $f : I \to I$ is a continuous map.  If $f$ has a period-$m$ point and if $m \prec n$, then $f$ also has a period-$n$ point.

\item[\rm{(2)}]
For each positive integer $n$ there exists a continuous map $f : I \to I$ that has a period-$n$ point but that has no period-$m$ point for any $m$ with $m \prec n$.

\item[\rm{(3)}]
There exists a continuous map $f : I \to I$ that has a period-$2^i$ point for $i = 0, 1, 2, \ldots$ and has no periodic point of any other period.
\end{itemize}}

\noindent
{\it Proof.}
If $f$ has period-$j$ points with $j \ge 3$ and odd, then by Theorem 2 $f$ has period-$(j+2)$ points and by Theorem 3 $f$ has period-$(2 \cdot 3)$ points.  If $f$ has period-$(2 \cdot j)$ points with $j \ge 3$ and odd, then $f^2$ has period-$j$ points.  According to Theorem 2, $f^2$ has period-$(j+2)$ points, which implies that $f$ has either period-$(j+2)$ points or period-$(2 \cdot (j + 2))$ points.  If $f$ has period-$(j + 2)$ points, then in view of Theorem 3 $f$ also has period-$(2 \cdot (j+2))$ points.  In either case, $f$ has period-$(2 \cdot (j+2))$ points.  On the other hand, since $f^2$ has period-$j$ points, by Theorem 3 $f^2$ has period-$(2 \cdot 3)$ points, hence $f$ has period-$(2^2 \cdot 3)$ points.  Now if $f$ has period-$(2^k \cdot j)$ points with $j \ge 3$ and odd and if $k \ge 2$, then $f^{2^{k-1}}$ has period-$(2 \cdot j)$ points.  It follows from what we have just proved that $f^{2^{k-1}}$ has period-$(2 \cdot (j+2))$ points and period-$(2^2 \cdot 3)$ points.  It follows that $f$ has period-$(2^k \cdot (j+2))$ points and period-$(2^{k+1} \cdot 3)$ points.  Furthermore, if $f$ has period-$(2^i \cdot j)$ points with $j \ge 3$ and odd and if $i \ge 0$, then $f^{2^i}$ has period-$j$ points.  When $\ell \ge i$, it is clear that $f^{2^\ell} = (f^{2^i})^{2^{\ell-i}}$ has period-$j$ points.  By Lemma 1, $f^{2^\ell}$ has period-2 points.  This implies that $f$ has period-$2^{\ell+1}$ points whenever $\ell \ge i$.  Finally, if $f$ has period-$2^k$ points for some $k \ge 2$, then $f^{2^{k-2}}$ has period-4 points.  Again by Lemma 1, $f^{2^{k-2}}$ has period-2 points, ensuring that $f$ has period-$2^{k-1}$ points.  This proves (1). 

For (2) and (3), it suffices to assume that $I = [0, 1]$.  Let $T(x) = 1 - |2x - 1|$ be the tent map on $I$.  Then for each positive integer $n$ the equation $T^n(x) = x$ has exactly $2^n$ distinct solutions in $I$.  It follows that $T$ has finitely many period-$n$ orbits.  Among these period-$n$ orbits, let $P_n$ be one with the smallest diameter $\max P_n - \min P_n$ (or one with the smallest $\max P_n$ or one with the largest $\min P_n$).  For any $x$ in $I$ let $T_n(x) = \min P_n$ if $T(x) \le \min P_n$, $T_n(x) = \max P_n$ if $T(x) \ge \max P_n$, and $T_n(x) = T(x)$ if $\min P_n \le T(x) \le \max P_n$.  It is then easy to see that $T_n$ has exactly one period-$n$ orbit (i.e., $P_n$) but has no period-$m$ orbit for any $m$ with $m \prec n$ in the Sharkovsky ordering.  Now let $Q_3$ be any period-3 orbit of $T$ of minimal diameter.  Then $[\min Q_3, \max Q_3]$ contains finitely many period-6 orbits of $T$.  If $Q_6$ is one of smallest diameter, then $[\min Q_6, \max Q_6]$ contains finitely many period-12 orbits of $T$.  We choose one, say $Q_{12}$, of minimal diameter and continue the process inductively.  Let $q_0 = \sup \{\min Q_{2^i \cdot 3} : i \ge 0 \}$ and $q_1 = \inf \{ \max Q_{2^i \cdot 3} : i \ge 0 \}$.  Let $T_{\infty}(x) = q_0$ if $T(x) \le q_0$, $T_{\infty}(x) = q_1$ if $T(x) \ge q_1$, and $T_{\infty}(x) = T(x)$ if $q_0 \le T(x) \le q_1$.  Then it is easy to check that $T_{\infty}$ has a period-$2^i$ point for $i = 0, 1, 2, \ldots$ but has no periodic point of any other period.  This establishes (2) and (3).
\hfill\sq

\noindent
{\bf Remark.} Lemma 1, although simple, has many applications.  In addition to Theorems 2 and 3 and the one stated in $\bold{[2}, {\rm sec.} 2$], we present two more such applications.  Assume that there exist a point $x_0$ in $I$ and an integer $n$ $(\ge 2)$ such that $f^n(x_0) \le x_0 < f(x_0)$.  Let $X = \{ f^k(x_0) : 0 \le k \le n-1 \}$, let $a = \max \{ x \in X : x_0 \le x$, and $f(x) > x \}$, and let $b$ be any point in $\{ x \in X : a < x \le f(a) \}$ such that $f(b) < a$.  Then $f(b) < a < b \le f(a)$.  Lemma 1 is therefore applicable to the situation.  We urge the interested reader to use Lemma 1 (with arguments similar to those used in the proofs of Theorems 2 and 3) to establish two facts.  First, if $f^n(x_0) \le x_0 < f(x_0)$ for some {\it odd} integer $n$ greater than one, then $f$ has period-$n$ points {\bf{\cite{li}}}.  Second, if $f$ is transitive (i.e., $f$ has a dense orbit) {\bf{\cite{hs}}} or, more generally, if the $\omega$-limit set $\omega(f,c)$ of some point $c$ ($x$ is in $\omega(f,c)$ if and only if $\lim_{k \to \infty} f^{n_k}(c) = x$ for some sequence of positive integers $n_k \to \infty$) contains both a fixed point $z$ and a point different from $z$, then $f$ has periodic points of all even periods and $f^2$ is turbulent (i.e., there exist closed subintervals $J_0$ and $J_1$ of $I$ with at most one point in common such that $f^2(J_0) \cap f^2(J_1) \supp J_0 \cup J_1$).

\noindent
{\bf ACKNOWLEDGMENTS.}  I would like to thank A. N. Sharkovsky and the referee for many suggestions that improved this note.

\end{document}